\documentclass{amsart}
\usepackage{amssymb}
\usepackage{enumerate}
\usepackage{url}

\theoremstyle{plain}
\newtheorem{prop}{Proposition}[section]
\newtheorem{theorem}[prop]{Theorem}

\newtheorem{lemma}[prop]{Lemma}
\newtheorem{conjecture}[prop]{Conjecture}
\theoremstyle{definition}
\newtheorem{definition}[prop]{Definition}
\newtheorem{remark}[prop]{Remark}
\newtheorem{example}[prop]{Example}

\newcommand{\N}{\mathbb N}
\newcommand{\R}{\mathbb R}
\newcommand{\axy}{\langle X,Y\rangle} 
\newcommand{\s}[1]{S_{#1}(X^2,Y^2)}
\newcommand{\csim}{\stackrel{\mathrm{cyc}}{\thicksim}}
\DeclareMathOperator{\tr}{tr}
\DeclareMathOperator{\ord}{ord}
\newcommand{\bm}{\mathbf}

\begin{document}
\title{Sums of hermitian squares as an approach\\ to the BMV conjecture}
\author{Sabine Burgdorf}
\address{Institut de Recherche Math\'ematique de Rennes,
         Universit\'e de Rennes I,
         Campus de Beaulieu,
         35042 Rennes cedex,
         France}
\email{sabine.burgdorf@univ-rennes1.fr}
\date{May 2009}

\begin{abstract}
Lieb and Seiringer stated in their reformulation of the Bessis-Moussa-Villani conjecture that all coefficients of the polynomial $p(t)=\tr[(\bm{\rm A}+\bm{\rm B})^m]$ are nonnegative whenever $\bm{\rm  A}$ and $\bm{\rm B}$ are any two positive semidefinite matrices of the same size. We will show that for all $m\in\N$ the coefficient of $t^4$ in $p(t)$ is nonnegative, using a connection to sums of Hermitian squares of non-commutative polynomials which has been established by Klep and Schweighofer. This implies by a famous result of Hillar that the coefficients of $t^k$ are nonnegative for $0\leq k\leq4$.
\end{abstract}

\maketitle

\section{Introduction}
The Bessis-Moussa-Villani (BMV) conjecture, originally stated as a problem of quantum statistical mechanics, has a 30 year long history. Since its introduction in 1975 \cite{bmv} many partial results have been given, see e.g. \cite{mou} for a review until 2000.  
The following reformulation of Lieb and Seiringer \cite{ls} is more capable to algebraic methods than the original one.
\begin{conjecture}[\bf(Bessis, Moussa, Villani)]\label{bmv}
For all positive semidefinite matrices $\bm{\rm A}$ and $\bm{\rm B}$ and all $m\in\N$, the polynomial $p(t):=\tr((\bm{\rm A}+t\bm{\rm B})^m)\in\R[t]$ has only nonnegative coefficients.
\end{conjecture}

The coefficient of $t^k$ in $p(t)$ for a given $m$ is the trace of $S_{m,k}(\bm{\rm A},\bm{\rm B})$, where $S_{m,k}(\bm{\rm A},\bm{\rm B})$ is the sum of all words of length $m$ in the letters $\bm{\rm A}$ and $\bm{\rm B}$ in which $\bm{\rm B}$ appears exactly $k$ times. For example $S_{4,2}(\bm{\rm A},\bm{\rm B})=\bm{\rm A}^2\bm{\rm B}^2+\bm{\rm A}\bm{\rm B}\bm{\rm A}\bm{\rm B}+\bm{\rm A}\bm{\rm B}^2\bm{\rm A}+\bm{\rm B}\bm{\rm A}\bm{\rm B}\bm{\rm A}+\bm{\rm B}^2\bm{\rm A}^2+\bm{\rm B}\bm{\rm A}^2\bm{\rm B}$.

In \cite{bmv} it has already been shown that  the BMV conjecture is true for $2\times 2$ matrices. Since for $0\leq k\leq 2$ or $m-2\leq k\leq m$ each word in $S_{m,k}(\bm{\rm A},\bm{\rm B})$ has nonnegative trace, as is easily seen, the conjecture is true for $m\leq 5$. Hillar and Johnson \cite{hj} verified the first nontrivial case $m=6, k=3$ for positive semidefinite $3\times 3$ matrices.
H\"agele \cite{hag} verified $m= 7$  which leads by a result of Hillar \cite{hil} to $m\leq 7$.  Further, Klep and Schweighofer \cite{ksbmv} derived that Conjecture \ref{bmv} is true for $m\leq 13$. Whereas all these results fix $m$ and consider arbitrary $k\leq m$, we take the opposite viewpoint, fix $k=4$ and let $m\in\N$ be arbitrary. We will give a proof that $$\tr(S_{m,4}(\bm{\rm A},\bm{\rm B}))\geq0$$ with no restrictions on $m$ or the matrix size of $\bm{\rm A}$ and $\bm{\rm B}$. 

A result of Hillar \cite{hil} then implies that it is true for all $k\leq 4$ and arbitrary $m$, in particular if $k=3$ which can't be shown directly by our method.

Using analytical methods, Fleischhack and Friedland \cite{ff} showed that for fixed  positive semidefinite $\bm{\rm A}, \bm{\rm B}$ and fixed $k$, the trace of  $S_{m,k}(\bm{\rm A},\bm{\rm B})$ is nonnegative whenever $m$ is big enough. Unfortunately, their lower bound on $m$ is dependent of $\bm{\rm A}$ and $\bm{\rm B}$. Otherwise this would imply the BMV conjecture.
 
To verify Conjecture \ref{bmv} it is sufficient to show the nonnegativity of $\tr(S_{m,k}(\bm{\rm A},\bm{\rm B}))$ for any two  positive semidefinite \textit{real} matrices $\bm{\rm A}$ and $\bm{\rm B}$ of the same size. Since further 
every positive semidefinite real matrix $\bm{\rm A}$ is decomposable as $\bm{\rm A}=\bm{\rm C}^2$ for some real matrix $\bm{\rm C}$ we specify our examination to $S_{m,k}(\bm{\rm C}^2,\bm{\rm D}^2)$ where $\bm{\rm C}$ and $\bm{\rm D}$ are any two real symmetric matrices of the same size. To work in an algebraic context we identify  $S_{m,k}(\bm{\rm C}^2,\bm{\rm D}^2)$ as a polynomial $S_{m,k}(X^2,Y^2)$ in two non-commuting variables $X$ and $Y$. 

For this let $\R\axy$ denote the unital associative $\R$-algebra freely generated by $X$ and $Y$. The elements of $\R\axy$
are polynomials in the non-commuting variables $X$ and $Y$ with coefficients in $\R$. An element $w$ of the monoid $\axy$, freely generated by $X,Y$, is called a \textit{word} and $w_{(i)}$ its $i$-th letter. An element of the form $aw$, where $0\neq a\in\R$ and $w\in\axy$, is called a \textit{monomial} and $a$ its \textit{coefficient}. We endow $\R\axy$ with the involution $p\mapsto p^*$
fixing $\R\cup\{X,Y\}$ pointwise. In particular, for each word $w\in\axy$, $w^*$ is its reverse. If $w^*=w$, $w$ is called a \textit{palindrome}. An element of the form $g^*g$ for some $g\in\R\axy$ is called a \textit{hermitian square}.

Using this terminology we define the polynomial $S_{m,k}(X,Y)$ as the polynomial in the variables $X$ and $Y$ as the sum of all monic monomials of total degree $m$ and degree $k$ in $Y$. Replacing $X$ and $Y$ by $X^2$ and $Y^2$ leads to the desired polynomial $\s{m,k}$, which results in $S_{m,k}(\bm{\rm A},\bm{\rm B})$ when we evaluate at symmetric matrices $\bm{\rm C}$ and $\bm{\rm D}$, satisfying $\bm{\rm C}^2=\bm{\rm A}$ and $\bm{\rm D}^2=\bm{\rm B}$.   

The invariance of the trace under cyclic permutations motivates the definition of cyclic equivalence \cite{ksbmv}. A cyclic permutation of a word $v$ of length $m$ is a map $\sigma$, where $\sigma(v)=v_{(\sigma(1))}v_{(\sigma(2))}\cdots v_{(\sigma(m))}$, for which there exists some $k\in\N$ such that $\sigma(i)=i+k \mod m$ for all $i=1,\dots,m$. For example $v_{(1)}v_{(2)}v_{(3)}\mapsto v_{(3)}v_{(1)}v_{(2)}$ is a cyclic permutation whereas $v_{(1)}v_{(2)}v_{(3)}\mapsto v_{(3)}v_{(2)}v_{(1)}$ is not.

\begin{definition}
Two words $v,w\in\axy$ are called \textit{cyclically equivalent} ($v\csim w$) if $\sigma(v)=w$ for some cyclic permutation $\sigma$ of $v$. 

Two polynomials $f=\sum_{w}a_ww$ and $g=\sum_{w}b_ww$ with $a_w,b_w\in\R$ are cyclically equivalent if for each $v\in\axy$ the sums of coefficients of all words $w\in\axy$ which are cyclically equivalent to $v$ are equal, i.e.,  $\sum_{{w\csim v}}a_w=\sum_{{w\csim v}}b_w.$ This is equivalent to $f-g$ being a sum of commutators in $\R\axy$, where the commutator $[p,q]$ is defined by $[p,q]:=pq-qp$. 
\end{definition}

The polynomials $f=X^2YX+YX^3+2 X^2Y^2$ and $g=2 YX^3+2 YX^2Y$ are cyclically equivalent since $f-g=[X^2,YX]+[2 X^2Y,Y]$. Alternatively, the condition on the coefficients is easily checked as well.

\begin{definition}\label{cycsf} 
The \textit{order} ($\ord w$) of a word $w=w_{(1)}\cdots w_{(m)}$ of length $m$ is the smallest positive integer $k$, such that $w_{(i+k)}=w_{(i)} $ for all $i=1,\dots,m$ where we identify $w_{(i+k)}$ with $w_{(i+k-m)}$ if $i+k>m$. Thus cyclically equivalent words have the same order. It can also be defined as the smallest integer $k\geq1$ such that there exists a subword $v=v_{(1)}\cdots v_{(k)}$ of length $k$ with $w=v\cdots v=v^{m/k}$. The equivalence of these two definitions follows easily by induction over the length of the subword $v$.
\end{definition}

\begin{remark}\label{rts}
One obtains that the order of a word $w=v^{m/\ord(w)}$ in $\s{m,4}$ divides $m$. Further, since $Y^2$ appears the same number of times in every subword $v$, we get that $\frac{m}{\ord(w)}$ divides $4$. In particular $\ord(w)\in\{m,\frac m2, \frac m4\}\cap \N$.
\end{remark}

Our main result is the following. 
\begin{theorem}\label{shs}
For $k=4$ and $m\in\N$ the polynomial $\s{m,4}$ is cyclically equivalent to a sum of Hermitian squares. 
\end{theorem}

\begin{remark}\label{nonneg}
\begin{enumerate}[(i)]
\item H\"agele \cite{hag} has shown that $\s{6,3}$ cannot be cyclically equivalent to a sum of Hermitian squares of a certain special form. Landweber and Speer generalized this result to $k=3$ and $m\geq 6$ but $m\neq 11$ \cite{las}. Using this result a fact of Klep and Schweighofer \cite[Prop. 3.1]{ksbmv} shows that $\s{m,3}$ cannot be cyclically equivalent to any sum of Hermitian squares if $m\geq 6$ and $m\neq 11$. Therefore we are interested in the case $k=4$ and arbitrary $m\in\N$. 
\item Since a sum of Hermitian squares is positive semidefinite on all real symmetric matrices, Theorem \ref{shs} implies that $\tr(S_{m,4}(A,B))$, the coefficient of $t^4$ in $p(t)$, for all $m\in\N$ is nonnegative for all positive semidefinite matrices $A,B$. 
\end{enumerate}
\end{remark}

In the sequel we will present a proof of Theorem \ref{shs} by constructing a sum of Hermitian squares which is cyclically equivalent to $\s{m,4}$. By Remark \ref{rts} the order of words in $\s{m,4}$ divides $m$ and $4$. Thus, if $m$ is odd all words in $\s{m,4}$ have order $m$, whereas in the even case order $\frac m2$ and $\frac m4$ are also possible. Therefore we split the proof in two parts, $m$ odd and even, starting with the easier part, where $m$ is odd.

\section{Case $m$ odd}
To verify Theorem \ref{shs} it suffices to construct a sum of Hermitian squares $f$ which is cyclically equivalent to $\s{m,4}$. Let $m$ be fixed. Since $\s{m,4}$ is homogeneous in $X$ and $Y$, one can reduce the set of words in a decomposition as sum of Hermitian squares, as in the commutative case, to the set of words of half the degree in $X$ and $Y$. Thus we set
$$V=\{v\in\langle X,Y\rangle \mid \deg_{X}v=m-4, \deg_{Y}v=4 \}.$$
Further we define the subsets 
\begin{align*}
V_0&=\{v\in \{X^2,Y^2\}^{\frac{m-1}{2}}X \mid v=X^{k}Y^2X^{\ell}Y^2X^{k'+1}, k\leq k'\}\cap V ,\\ 
V_1&=\{v\in X\{X^2,Y^2\}^{\frac{m-1}{2}} \mid v=X^{k+1}Y^2X^{\ell}Y^2X^{k'}, k+1\leq k' \}\cap V.
\end{align*} 

We denote the possible exponents of $X$ in a word $v_i$ by $k_i,\ell_i$ and $k_i'$ such that for example
every $v_i\in V_0$ is of the form $v_i=X^{k_i}Y^2X^{\ell_i}Y^2X^{k_i'+1}$ and satisfies the condition 
$k_i+\ell_i+k_i'=m-5$ where $\ell_i,k_i,k_i'\in 2\N$ and $k_i\leq k_i'.$

The exponent $k_i$ (respectively $k_i+1$ if $v_i\in V_1$) is bounded by d, the highest possible even (respectively odd) number which is less than or equal to $\frac{m-5}2$, thus the maximum of these bounds is in any case $\frac{m-5}2$. \\

Now, we will construct a sum of Hermitian squares $f$. For given $k\in\N$ let $k(2)$ denote the remainder of $k$ modulo 2. Then we group the words $v_i\in V_0$ (respectively $V_1$) according to $k_i$. For every $k=0,1,2,\dots,\frac{m-5}2$ we add all words $v_i\in V_{k(2)}$ with $k_i+k(2)=k$ and obtain a polynomial $f_{k}$. By construction
all words in $f_k^*f_k$ have even exponents in $X$ and $Y$. Finally, we set
\begin{equation}
   f:=m \sum_{k=0}^{\frac{m-5}2}f_{k}^\ast f_{k}.
\end{equation}

\begin{example}
\begin{enumerate}[(a)]
\item $m=7$: We have $V_0=\{Y^2X^2Y^2X,Y^4X^3\}$ and $V_1=\{XY^4X^2\}$ which leads to 
\begin{align*}
f_0= Y^2X^2Y^2X+Y^4X^3 \quad\text{and}\quad f_1= XY^4X^2 
\end{align*}
and finally
\begin{align*}
\s{7,4}&\csim 7\left(f_0^* f_0+f_1^* f_1\right) \\
 &=7(XY^2X^2Y^4X^2Y^2X+ XY^2X^2Y^6X^3+X^3Y^6X^2Y^2X+ X^3Y^8X^3\\
 &\;\qquad + XY^4X^4Y^4X).
\end{align*} 
This representation is of the same kind as the one given by H\"agele in \cite{hag}.

\item $m=9$: Since $V_0=\{Y^2X^2Y^2X^3,Y^4X^5,X^2Y^4X^3,Y^2X^4Y^2X\}$ and 
$V_1=\{XY^4X^4,XY^2X^2Y^2X^2 \}$ we get by construction
\begin{align*}
f_0&= Y^2X^2Y^2X^3+Y^4X^5+Y^2X^4Y^2X, \\
f_1&= XY^2X^2Y^2X^2+XY^4X^4 \quad\text{and}\\
f_2&= X^2Y^4X^3. 
\end{align*}
One easily checks $\s{9,4}\csim 9\left(f_0^* f_0+f_1^* f_1+f_2^* f_2\right).$ 
\end{enumerate}
\end{example}

We will prove that $f$ is the desired sum of Hermitian squares in two steps. First all words appearing in $f$ will be shown to be pairwise cyclically inequivalent. By construction each word in $f$ appears in $\s{m,4}$ and has order $m$. Since up to cyclic equivalence each word in $\s{m,4}$ appears $m$ times, it suffices to show that the sums of coefficients in both polynomials are the same. 

\begin{remark}\label{geom}
To compare two words appearing in $f$ with respect to cyclic equivalence we use the following method. Since $Y^2$ appears exactly four times in each word $w$ of $f$, we know $w=X^{n_0}Y^2X^{n_1}Y^2X^{n_2}Y^2X^{n_3}Y^2X^{n_4'}$ for some $n_0,\dots,n_3,n_4'$. Further $w$ is cyclically equivalent to $\tilde w:=Y^2X^{n_1}Y^2X^{n_2}Y^2X^{n_3}Y^2X^{n_4}$ where $n_4=n_4'+n_0$, i.e. $\tilde w$ consists of four groups $Y^2X^{n_i}$. Let $w'$ be another word with exponents $m_i$, i.e., $\tilde w':=Y^2X^{m_1}Y^2X^{m_2}Y^2X^{m_3}Y^2X^{m_4}$. Then $\tilde w$ and $\tilde w'$ are the same or $n_i=m_{i-j}\;(i-j \mod 4)$ for $i=1,\dots,4$ and $j=1,2,3$, which can be obtained by ``rotating'' $\tilde w'$ $j$ times, i.e., for $j=1$ one shifts the first group $Y^2X^{m_1}$ to the end, for $j=2$ one shifts also the second group to the end and so on, thus $m_i$ becomes $m_{i-j} $.
 
For simplicity we use the fact that rotating three times is the same as rotating once in the reverse direction, i.e., shifting the group $Y^2X^{m_4}$ to the beginning. Thus rotating $w'$ three times is the same as fixing $w'$ and rotating $w$ once. Therefore we can omit $j=3$ by symmetry.
\end{remark}

\begin{lemma}\label{ncsim}
All words appearing in $f$ are pairwise cyclically inequivalent. 
\end{lemma}

\begin{proof}
By construction a word $w$ in $f$ is either a word in $\sum_{2k} f_{2k}^*f_{2k}$ thus of the form $w=v_1^*v_2$ where $v_1,v_2\in V_0$ and $k_1=k_2$, i.e., $$w=X^{k_1'+1}Y^2X^{\ell_1}Y^2X^{2k_1}Y^2X^{\ell_2}Y^2X^{k_2'+1}\csim Y^2X^{\ell_1}Y^2X^{2k_1}Y^2X^{\ell_2}Y^2X^{k_1'+k_2'+2}.$$ Or it is a word in $\sum_{2k} f_{2k+1}^*f_{2k+1}$ thus of the form $w=v_1^*v_2$ where $v_1,v_2\in V_1$ and $k_1=k_2$. The same is true for any other word $w'=v_3^*v_4$.  As is easily seen $\tilde w=\tilde w'$ is only possible if $v_1=v_3$ and $v_2=v_4$. We are left with the following cases.

If $w$ and $w'$ are words in $\sum_{2k} f_{2k}^*f_{2k}$ which are cyclically equivalent then we have to consider
\begin{enumerate}[(a)]
\item $\ell_1=2k_3,  \quad 2k_1= \ell_4, \quad \ell_2=k_3'+k_4'+2,  \quad k_1'+k_2'+2= \ell_3$ or
\item $\ell_1=\ell_4,  \quad 2k_1=k_3'+k_4'+2,  \quad \ell_2= \ell_3,   \quad k_1'+k_2'+2= 2k_3$.
\end{enumerate}
In (a) $2k_3+k_1+k_1'=\ell_1+k_1+ k_1'=\ell_3+k_3+ k_3'=k_1'+k_2'+2+k_3+ k_3'$ leads to $k_1+k_3=k_2'+k_3'+2$ contradicting $k_1+k_3\leq k_2'+k_3'<k_2'+k_3'+2$.  Subcase (b) leads to $2k_1=k_3'+k_4'+2> 2 k_3=k_1'+k_2'+2>2k_1$, which is not possible. 

The case that $w,w'$ are words in $\sum_{2k} f_{2k+1}^*f_{2k+1}$ works the same way. 

If $w$ is a word in $\sum_{2k} f_{2k}^*f_{2k}$ and $w'$ a word in $\sum_{2k} f_{2k+1}^*f_{2k+1}$, then we have
\begin{enumerate}[(a)]
\item $ \ell_1=k_3'+k_4', \quad 2k_1= \ell_3, \quad \ell_2=2k_3+2, \quad k_1'+k_2'+2= \ell_4$ or
\item $ \ell_1=\ell_3, \quad 2k_1=2k_3'+2, \quad \ell_2= \ell_4, \quad k_1'+k_2'+2= k_3'+k_4'$.
\end{enumerate}
In (a) $k_3'+k_4'+k_1+k_1'=\ell_1+k_1+k_1'=\ell_3+k_3+k_3'=2k_1+k_3+k_3'$ leads to $k_1'+k_4'=k_1+k_3=k_1+k_4<k_1'+k_4'$. Subcase (b) contradicts $k_1,k_3'\in 2\N$. 

If $w$ is a word in $\sum_{2k} f_{2k+1}^*f_{2k+1}$ and $w'$ in $\sum_{2k} f_{2k}^*f_{2k}$, we exchange $w$ and $w'$.\\

Summarizing, despite the trivial case that $w$ and $w'$ are constructed by the same subwords $v_i$, they cannot be cyclically equivalent.
\end{proof}

Thus every word in $f$ has its order $m$ as coefficient. Since up to cyclic equivalence this is the same in $\s{m,4}$, we are done by the following lemma. 

\begin{lemma}\label{anzahl}
The number of pairwise cyclically inequivalent words in $f$ is the same as in $\s{m,4}$.
\end{lemma}

\begin{proof}
$\s{m,4}$ contains $\binom{m}{4}$ words. Since each word has order $m$, there are $$\frac 1m\binom{m}{4}=\frac 1{6} (\frac{m-3}2)(\frac{m-1}2)(m-2)$$ pairwise cyclically inequivalent words in $\s{m,4}$.

Let $k\in\N$ be fixed. Then $f_{k}$ consists of $\frac{m-3}2-k$ different words. For example, if $k$ is even then there are $\frac12 (m-5-k_1)+1$ possibilities for $k_1,\ell_1,k_1'\in2\N$ with $\ell_1+k_1'=m-5-k_1$ (namely $k_1'=m-5-k_1-\ell_1, \ell_1=0,2,\dots m-5-k_1$), the restriction $k_1\leq k_1'$ of $V_0$ excludes $\frac{k_1}2$ possibilities.

Thus the number of words in $f$ is given by
\begin{align*} 
&\sum_{k=0}^{\frac{m-5}2}(\frac{m-3}2-k)^2 =\sum_{k=0}^{\frac{m-3}2} k^2 =\frac16(\frac{m-3}2)(\frac{m-1}2)(m-2).
\end{align*} 
\end{proof}

\begin{remark}\label{v3odd}
 After we had finished the proof of this case, we heard of the recent work of Landweber and Speer \cite{las} who proved the same result (for odd $m$) by quite similar techniques; but they haven't investigated the case where $m$ is even. They found a sum of Hermitian squares which only consists of words $w$ in 
$$V_2:=\{v\in X\{X^2,Y^2\}^{\frac{m-1}{2}} \mid v=X^{k+1}Y^2X^{l}Y^2X^{k'} \}\cap V.$$ 
Let $v_i=w_i X\in V_0$ for $i=1,2$. Starting with $f$ and using 
$$v_1^\ast v_2=(Xw_1)^*Xw_2= w_1^*XXw_2\csim Xw_2w_1^*X=(w_2^*X)^*(w_1^*X)=(v_2^\ast)^\ast (v_1^\ast)$$ and $V_1\subseteq V_2$ 
leads to a sum of Hermitian squares $\tilde f$ which is exactly the representation found by Landweber and Speer. 
 
This result agrees with the more general Proposition 3.1 in \cite{ksbmv} which in particular states that independent of $k$ in the case $m$ odd once one has found a representation as sum of Hermitian squares one can also find a representation using only of words of $V_2$. 
\end{remark}

\section{Case $m$ even}

Since words in $\s{m,4}$ now have order $m,\frac m2$ or $\frac m4$, the constructed polynomial $f$  of the last section is further not cyclically equivalent to $\s{m,4}$. Thus we will add weights on the words in our construction to respect the different orders.  \\

Let $m$ be fixed and $V=\{v\in\langle X,Y\rangle|\deg_{X}v=m-4, \deg_{Y}v=4 \}.$
Further we define the subsets 
\begin{align*}
V_0&=\{v\in\{X^2,Y^2\}^{\frac{m}{2}}\mid v=X^{k}Y^2X^{\ell}Y^2X^{k'}, k\leq k' \}\cap V,\\
V_1&=\{v\in X\{X^2,Y^2\}^{\frac{m-2}{2}}X\mid v=X^{k+1}Y^2X^{\ell}Y^2X^{k'+1}, k\leq k' \}\cap V. 
\end{align*}

To distinguish even and odd exponents, we define $\hat{k_i}:=k_i+1$ and $\hat{k_i'}:=k_i'+1$. Then every $v_i\in V_0$ is of the form $v_i=X^{k_i}Y^2X^{\ell_i}Y^2X^{k_i'+1}$ and satisfies $k_i+\ell_i+k_i'=m-4$ where $\ell_i,k_i,k_i'\in 2\N$ and  $k_i\leq k_i'$, 
whereas every $v_i\in V_1$ satisfies $\hat k_i+\ell_i+\hat k_i'=m-4$. Thus the maximal possible exponent $k_i$ respectively $\hat k_i$ (if $m$ is not divisible by 4) is given by $\frac{m-4}2$.\\  

Now we construct our desired sum of Hermitian squares as follows. 
Let $k\in\N$ and let $k(2)$ denote the remainder of $k$ modulo $2$. For every $k=0,1,2,\dots\frac{m-4}2$ we add all words $v_i\in V_{k(2)}$ with $k_i+k(2)=k$ as in the case where $m$ is odd, but we weight the words with $k_i<k_i'$ with coefficient $1$ and the words with $k_i=k_i'$ with coefficient $\frac12$. This leads to a polynomial $f_{k}$ which contains exactly one word with coefficient $\frac12$ whereas all other coefficients are $1$. Finally we set 
\begin{equation}
   f:=m \sum_{k=0}^{\frac{m-4}2}f_{k}^\ast f_{k}.
\end{equation}
\begin{example}$m=8$: 
We have $V_0=\{Y^2X^2Y^2X^2,Y^4X^4,X^2Y^4X^2,Y^2X^4Y^2\}$
and $V_1=\{XY^4X^3,XY^2X^2Y^2X\}$ which leads to 
\begin{align*}
f_0&= Y^2X^2Y^2X^2+ Y^4X^4+\frac12 Y^2X^4Y^2 \\
f_1&= XY^4X^3+\frac12 XY^2X^2Y^2X \quad\text{and}\\
f_2&=\frac12 X^2Y^4X^2. 
\end{align*}
Then one easily verifies
$\s{8,4}\csim 8(f_0^\ast f_0+f_1^\ast f_1+f_2^\ast f_2)$.

Now for example the words $w=(Y^2X^2Y^2X^2)^*(Y^2X^4Y^2)$ in $f_0^*f_0$ and $w'=(XY^2X^2Y^2X)^*(XY^4X^3)$ in $f_1^*f_1$ are cyclically equivalent but due to our weights their coefficients sum up to $\ord(w)=8$.
\end{example}

The proof of cyclic equivalence works similarly as in the case where $m$ is odd. But since there are now cyclically equivalent words appearing in $f$, we have to calculate more carefully. We will show first that the sum of coefficients of cyclically equivalent words in $f$ is less than or equal to their order. Since each word in $f$ appears in $\s{m,4}$ we will finish by showing that the sums of coefficients are equal in both representations.
\begin{lemma}\label{coeff}
The sum of coefficients of cyclically equivalent words in $f$ is less than or equal to the order of the corresponding words. 
\end{lemma}

\begin{proof}
We will use the same method as explained in Remark \ref{geom} of the last section. Let $w,w'$ be two different words appearing in $f$ and $w\csim w'$. 

If $w$ and $w'$ are words in $\sum f_{2k}^*f_{2k}$.
Then either $w$ and $w'$ are equal or one of the following subcases holds:
\begin{enumerate}[(a)]
\item $ \ell_1=2k_3, \quad  2k_1= \ell_4, \quad  \ell_2= k_3'+k_4',  \quad k_1'+k_2'= \ell_3$
\item $ \ell_1=\ell_4,  \quad 2k_1= k_3'+k_4',  \quad \ell_2=\ell_3, \quad  k_1'+k_2'= 2k_3$
\end{enumerate}

In subcase (a) we obtain $k_3+k_1=k_2'+k_3'$ from $2 k_3+k_1+k_1'=\ell_1+k_1+k_1'=\ell_3+k_3+k_3'=k_1'+k_2'+k_3+k_3'$, thus $k_1=k_2'$ and $k_3=k_3'$. Further we obtain from $\ell_1+k_1+k_1'=\ell_4+k_3+k_4'$ that $k_1'-k_1=k_4'-k_3$. 
In (b) we obtain $2k_1=k_3'+k_4'\geq 2k_3=k_1'+k_2'\geq 2k_1$, thus equality holds, which leads to $w=w'$.

The other cases work in the same way by replacing $k_i$ by $\hat k_i$ whenever $w$ is a word in $\sum f_{2k+1}^*f_{2k+1}$ and $k_i'$ respectively if $w'$ is a word in $\sum f_{2k+1}^*f_{2k+1}$. If $w$ and $w'$ are not in the same set $\sum f_{2k}^*f_{2k}$ or $\sum f_{2k+1}^*f_{2k+1}$, then they obviously cannot be equal.   

Summarizing, we derife that when $w\csim w'$ but $w\neq w'$ then 
$k_1=k_2',\; k_3=k_3'$ and $k_1'-k_1= k_4'-k_3$ or by symmetry (confer Remark \ref{geom})
$k_3=k_4',\; k_1=k_1'$ and $k_3'-k_3=k_2'-k_1$ 
holds, where the first set of equations describes the words which differ by one rotation, and the second set describes the case of three rotations.

Assuming, there are two different words $w',w''$ both cyclically equivalent to $w$. Then all three are pairwise cyclically equivalent and at least two of them (for example $w',w''$) are in  $\sum_k f_{k}^*f_{k}$ ($k$ even or odd). Thus each of them satisfies one set of equations, but then $w'$ and $w''$ differ by two rotations, which leads to equality (subcase (b)). Therefore there are at most two words in $f$ which are pairwise cyclically equivalent. 

To conclude the proof, if  $w=v_1^*v_2$ with $k_1=k_1'=k_2=\frac{m-4}4'$ then $\ell_1=m-4-2k_1=\frac{m-4}2=\ell_2$, thus $w$ has order $\frac m4$ which is equal to the coefficient of $w$ in $f$. A cyclically equivalent word $w'=v_3^*v_4$ has to satisfy $k_3=k_3'=k_4'$ and $2k_3=\ell_1=2k_1$ which leads to $w=w'$. Therefore there is no other word $w'\csim w$ in $f$. 
In all other cases the coefficient of $w$ is half of the order of $w$. Since there are at most two pairwise cyclically equivalent words we are done.   
\end{proof}

\begin{lemma}\label{anzahl2}
The sum of coefficients in both polynomials is the same.
\end{lemma}

\begin{proof}
The sum of coefficients in $\s{m,4}$ is $\binom{m}{4}=\frac 1{24} m(m-1)(m-2)(m-3).$ 

For every $k=0,1,2\dots,\frac{m-4}2$ each polynomial $f_{k}$ has one word with coefficient $\frac{1}2$ and $\frac{m-4}2-k$ times coefficient $1$. Thus the sum of coefficients in $f$ is given by \vspace{-0.3cm}
\begin{align*} 
&m \sum_{k=0}^{\frac{m-4}2}\big(\frac{m-4}2-k+\frac 12\big)^2=\frac{m(m-2)}8+m\sum_{k=0}^{\frac{m-4}2}(k^2+k) \\
&=\frac {m}{24}\big(3(m-2)+(m-4)(m-2)m \big)
=\frac 1{24} m(m-1)(m-2)(m-3).
\end{align*} 
\end{proof}

\section{Concluding Remarks} 

\begin{enumerate}[(a)]
\item To get an idea how sums of Hermitian squares which are cyclically equivalent to $\s{m,4}$ might look like, we used numerical computations extending those done by Klep and Schweighofer \cite{ksbmv}. In particular we used NCAlgebra \cite{nca}, YALMIP \cite{yal} and SeDuMi \cite{sdm} as the starting point of our investigation.  

\item \label{v3even} As in the case $m$ odd one might consider $V_2=\{v\in \{X^2,Y^2\}^{\frac{m}{2}}\}\cap V $ if $m$ is even. Then one can find a much more complicated sum of Hermitian squares which is cyclically equivalent to $\s{m,4}$ and consists just of words in $V_2$ if $m(4)=2$, i.e., $m$ is even but not divisible by $4$. 

Since all words are words in the letters $X^2$ and $Y^2$ one obtains by substitution that $S_{m,4}(X,Y)$ is cyclically equivalent to a sum of Hermitian squares, which implies $\tr(S_{m,4}(\bm{\rm A},\bm{\rm B}))\geq 0$ for {\it all} real Hermitian matrices $\bm{\rm A},\bm{\rm B}$ of the same size. This result has recently, independently been found by Collins, Dykema and Torres-Ayala \cite{cdt}.

\item Landweber and Speer \cite{las} showed that despite a few exceptions (which all have been solved) one cannot find a representation of $\s{m,k}$ as a sum of Hermitian squares if $m$ or $k$ is odd.
But they have no negative results if $m$ and $k$ are both even. This gap has recently been filled by  Collins, Dykema and Torres-Ayala \cite{cdt} who proved that, despite the case (16,8), $S_{m,k}(X^2,Y^2)$ is not cyclically equivalent to a sum of Hermitian squares if $m-6\geq k\geq 6$ and $m\geq16$. Thus this approach cannot proof the BMV conjecture.  
\end{enumerate}

\subsection*{Acknowledgements}
I would like to thank M. Schweighofer for introducing me to the BMV-conjecture and the methods of semidefinite programming which he and I. Klep have used in \cite{ksbmv}, further I thank him for the helpful discussions on this topic. I also would like to thank I. Klep, P. Landweber and E. Speer for their helpful comments on previous versions.

\end{document}